\documentclass[10pt,a4paper]{amsart}
\usepackage[utf8]{inputenc}
\usepackage{amsmath}
\usepackage{amsfonts}
\usepackage{amssymb}
\usepackage{enumerate}
\usepackage{graphicx}
\usepackage{hyperref}
\usepackage{amsthm}
\usepackage{hhline}
\usepackage{mathdots}
\usepackage{xcolor}
\usepackage{colortbl}
\usepackage{listings}

\usepackage{mathtools}
\usepackage{subfigure}
\numberwithin{table}{section}
\usepackage{longtable}
\usepackage{tabularx}
\newcolumntype{C}[1]{>{\centering\arraybackslash}p{#1}}
\usepackage{booktabs}
\usepackage[english]{babel}

\calclayout
\newtheorem{prop}{Proposition}[section]
\newtheorem{theorem}[prop]{Theorem}
\newtheorem{lemma}[prop]{Lemma}
\theoremstyle{remark}
\newtheorem{rem}[prop]{Remark}
\theoremstyle{definition}

\newtheorem{condition}[prop]{Condition}
\DeclareMathOperator*{\Inn}{Inn}

\DeclareMathOperator*{\Aut}{Aut}

\DeclareMathOperator*{\SO}{SO}
\DeclareMathOperator*{\ord}{ord}
\DeclareMathOperator*{\Lin}{Lin}
\DeclareMathOperator{\Irr}{Irr}
\DeclareMathOperator{\Inf}{Inf}
\DeclareMathOperator{\SymM}{SymM}
\DeclareMathOperator{\Out}{Out}

\DeclareMathOperator*{\Gal}{Gal}
\DeclareMathOperator*{\ab}{ab}
\DeclareMathOperator*{\PSU}{PSU}
\DeclareMathOperator*{\OO}{O}

\DeclareMathOperator*{\GL}{GL}
\DeclareMathOperator*{\Sp}{Sp}

\DeclareMathOperator*{\PSL}{PSL}
\DeclareMathOperator*{\Res}{Res}
\DeclareMathOperator{\Ind}{Ind}

\DeclareMathOperator*{\Syl}{Syl}

\DeclareMathOperator*{\Def}{Def}
\DeclareMathOperator*{\pr}{pr}

\begin{document}
\title[On the inductive McKay--Navarro condition in the defining characteristic]{On the inductive McKay--Navarro condition for finite groups of Lie type in their defining characteristic}
\author{Birte Johansson}
\address{FB Mathematik, TU Kaiserslautern, Postfach 3049, 67653 Kaiserslautern, Germany.}
\email{johansson@mathematik.uni-kl.de}
\subjclass[2010]{20C33}
\keywords{groups of Lie type, local-global conjectures, McKay conjecture}
\date{\today}

\begin{abstract}
The McKay--Navarro conjecture is a refinement of the McKay conjecture that additionally takes the action of some Galois automorphisms into account. We verify the inductive McKay--Navarro condition in the defining characteristic for the finite groups of Lie type with exceptional graph automorphisms, the Suzuki and Ree groups, $\mathsf{B}_n(2)$ ($n \geq 2$), and the  groups of Lie type with non-generic Schur multiplier.
This completes the verification of the inductive McKay--Navarro condition for the finite groups of Lie type in their defining characteristic. 
\end{abstract}
\maketitle

\section{Introduction}
One of the main problems in the representation theory of finite groups is the study of so-called local-global conjectures. These conjectures assert a relation between the representation theory of a finite group $G$ and, for some prime $p$, of a local subgroup that is the normalizer of a $p$-subgroup of $G$. Many of these local-global conjectures are refinements of the McKay conjecture.

Let $G$ be a finite group, $p$ a prime, and $P$ a Sylow $p$-subgroup of $G$. The McKay conjecture claims that there is a bijection between the irreducible complex characters with degree not divisible by $p$ of $G$ and of the normalizer $N_G(P)$. Navarro refined the conjecture and proposed that this bijection can be chosen such that the same number of characters is fixed under the action of certain Galois automorphisms \cite{navarro2004mckayrefinement}. This is called the McKay--Navarro or Galois--McKay conjecture.

Navarro, Späth, and Vallejo reduced the McKay--Navarro conjecture to a problem about simple groups in \cite{navarro2019reduction}. If the inductive McKay--Navarro condition \cite[Definition 3.1]{navarro2019reduction} (originally called inductive Galois--McKay condition by the authors) is satisfied for all simple groups, then the McKay--Navarro conjecture itself holds for all groups.  In \cite{ruhstorfer2017navarro}, Ruhstorfer showed that the inductive McKay--Navarro condition is true for most groups of Lie type in their defining characteristic. 
We verify the inductive McKay--Navarro condition for the remaining families of groups of Lie type in their defining characteristic.
\begin{theorem}
The inductive McKay--Navarro condition \cite[Definition 3.1]{navarro2019reduction} is satisfied in the defining characteristic for the groups $\mathsf{B}_2(2^i)$, $\mathsf{G}_2(3^i)$, 
$\mathsf{F}_4(2^i)$, $\mathsf{B}_n(2)$ for integers $i \geq 1$, $n \geq 2$, the Suzuki and Ree groups, as well as for $\mathsf{G}_2(2)'$, $^2\mathsf{G}_2(3)'$, $^2\mathsf{F}_4(2)'$, and the simple groups of Lie type with non-generic Schur multiplier.
\end{theorem}
Together with \cite[Theorem 1.3]{ruhstorfer2017navarro}, this completes the verification of the inductive McKay--Navarro condition for finite groups of Lie type in their defining characteristic.

\subsection*{Acknowledgment} I would like to thank my PhD advisor Gunter Malle for his suggestions and comments on earlier versions as well as Lucas Ruhstorfer for his advice and providing me with an earlier updated version of his paper. Furthermore, I am thankful to the anonymous referee for helpful comments and suggesting a correction of the proof of Lemma \ref{lemmaHG}. This work was financially supported by the SFB-TRR 195 of the German Research Foundation (DFG).

\section{Preliminaries}
\subsection{Notation} 
Let $p$ be a prime and $G$ a finite group. We call a character of $G$ with degree not divisible by $p$ a $p'$-character and denote the set of complex irreducible $p'$-characters of $G$ by ${\Irr}_{p'}(G)$. The set of linear characters is denoted by $\Lin(G)$.

In the following, let $\mathcal{H}$ be the subgroup of $\mathcal{G}: =\Gal(\mathbb{Q}^{\ab}/\mathbb{Q})$ consisting of those $\sigma \in \mathcal{G}$ such that there exists an $n \in \mathbb{Z}$ with $\sigma(\zeta)=\zeta^{p^n}$  for every root of unity $\zeta$ of order not divisible by $p$. Here, the field $\mathbb{Q}^{\ab}$ denotes the subfield of $\mathbb{C}$ generated by all roots of unity.  
The group $\mathcal{H}$ acts naturally on characters and (projective) representations by applying the elements of $\mathcal{H}$ to the character values and matrix entries, respectively. 

We denote the automorphism group of $G$ by $\Aut(G)$, the inner automorphism group by $\Inn(G) \cong G/Z(G)$, and the outer automorphism group by $\Out(G)= \Aut(G)/\Inn(G)$. If it is clear from the context, we also use $\Inn(H)$ to denote the inner automorphisms of $G$ that arise from conjugation with elements of $H \leq G$.
Let $\psi$ be a character of a normal subgroup $N \unlhd G$ and  $\kappa \in \Aut(G)$. Then, we denote by $\psi^\kappa$ the character defined by $\psi^\kappa(x)=\psi(x^{\kappa^{-1}})$ for all $x \in N$, using exponential notation for the image of $x$ under the automorphism $\kappa$. For $g \in G$, we often use the element itself to denote the conjugation by $g$, i.e. $\psi^g(x)=\psi(x^{g^{-1}})=\psi(gxg^{-1})$ for all $x \in N$. The stabilizer of $\psi$ in $G$ under this conjugation action is denoted by $G_\psi$. In the same way, automorphisms of a group act on the (projective) representations of its normal subgroups. 
For a subgroup $A \subseteq \Aut(G)$, we denote the setwise stabilizer of the orbit $\psi^\mathcal{H}$ under the action of $A$ by $A_{\psi^\mathcal{H}}:=\{\kappa \in A \mid \kappa(\psi^\mathcal{H}) \subseteq \psi^\mathcal{H} \}$.

\subsection{Inductive McKay--Navarro condition} 
We need the notion of associated projective representations and a way to describe how they behave under conjugation and Galois automorphisms.
\begin{lemma}\cite[Corollary 1.2 and Lemma 1.4] {navarro2019reduction} \label{lemmamu_a}
Let $N \unlhd G$ and $\psi$ be a character of $N$. 
\begin{enumerate}[(a)]
\item There exists a projective representation $\mathcal{P}$ of $G_\psi$ over $\mathbb{Q}^{\ab}$ such that 
its restriction to $N$ affords $\psi$ and
$\mathcal{P}(ng)=\mathcal{P}(n)\mathcal{P}(g)$, $\mathcal{P}(gn)=\mathcal{P}(g)\mathcal{P}(n)$ for all $n \in N$ and $g \in G_\psi$. 
\item \label{defmu_a} Let $\mathcal{P}$ be such a projective representation and $(g, \sigma) \in G \times \Gal(\mathbb{Q}^{\ab}/\mathbb{Q})$ with $\psi^{g \sigma}:=(\psi^g)^\sigma=\psi$. Then there exists a unique function $\mu_{g  \sigma}: G_\psi \rightarrow \mathbb{C}^\times$ constant on cosets of $N$ with $\mu_{g \sigma}(1)=1$ such that $$\mathcal{P}^{g \sigma } := (\mathcal{P}^g)^\sigma \sim \mu_{g \sigma} \mathcal{P}$$ where $\sim$ denotes similarity between the projective representations. 
\end{enumerate}
\end{lemma}
We say that such a projective representation $\mathcal{P}$ is \textit{associated to $\psi$}. In the following, we use $\mu_{g  \sigma}$ to denote the transition function defined here. The lemma also allows us to assume that all occurring (projective) representations are defined over $\mathbb{Q}^{\ab}$.

The inductive McKay--Navarro condition from \cite[Definition 3.1]{navarro2019reduction} can be stated in the following way.
\begin{condition}\label{imck}
For a finite non-abelian simple group $S$ and $p$ a prime dividing $|S|$, let $G$ be a universal covering group of $S$, $R \in \Syl_p(G)$, and $\Gamma :=\Aut(G)_R :=\{f \in \Aut(G) \mid f(R) = R\}$. Then, $S$ satisfies the \textit{inductive McKay--Navarro condition for $p$} if the following holds for all $\psi \in {\Irr}_{p'}(G)$:
\begin{enumerate}[(01)]
\item[(1)] \textit{(Equivariant bijection)} There exists a $\Gamma$-stable subgroup $N_G(R) \subseteq N \subsetneq G$ and a $\Gamma \times \mathcal{H}$-equivariant bijection $$\Omega: {\Irr}_{p'}(G) \rightarrow {\Irr}_{p'}(N).$$ 
\item[(2A)] \textit{(Extension condition)} There exist projective representations $\mathcal{P}$ of $G\rtimes \Gamma_{\psi}$ and $\mathcal{P}'$ of $N \rtimes \Gamma_{\Omega(\psi)}$ with entries in $\mathbb{Q}^{\ab}$ associated to $\psi$ respectively $\Omega(\psi)$ such that the associated factor sets $\alpha, \alpha'$ take roots of unity values, coincide on $(N \rtimes \Gamma_{\psi}) \times (N \rtimes \Gamma_{\psi})$, and the scalar matrices $\mathcal{P}(c)$, $\mathcal{P}'(c)$ correspond to the same scalar for all $c \in C_{G\rtimes \Gamma_{\psi^\mathcal{H}}}(N)$.
\item[(2B)]
Further, $\mu_a$ and $\mu_a'$ agree on $N \rtimes \Gamma_{\psi}$ for all $a \in ((N \rtimes \Gamma_{\psi^\mathcal{H}}) \times \mathcal{H})_{\psi}$.
\end{enumerate}
\end{condition}
\begin{rem}
\begin{enumerate}[(a)] \label{rem1}
\item 
The extension condition ensures that we have $$(G \rtimes \Gamma_{\psi^\mathcal{H}}, G, \psi) \geq_c (N \rtimes \Gamma_{\psi^\mathcal{H}}, N, \Omega(\psi))$$ in the language of $\mathcal{H}$-triples as introduced in \cite[Definition 1.5]{navarro2019reduction}. 
\item \label{remmckaydefining}
Note that we already know that the inductive McKay condition from \cite[Section 10]{isaacsmallenavarro2007} is satisfied for groups of Lie type in their defining characteristic \cite[Theorem 1.1]{spaeth2012mckaydefining}. Thus, it remains to verify the $\mathcal{H}$-equivariance of the existing $\Gamma$-equivariant bijections $\Irr_{p'}(G) \rightarrow \Irr_{p'}(N)$ and show that $\Gamma \times \mathcal{H}$ acts on the associated projective representations in the same way for the local and global characters.
\item \label{remgallagher}
Let $\mathcal{X}$ be a representation affording $\psi \in {\Irr}_{p'}(G)$. We can extend it canonically to the inner automorphisms in $\Gamma$ by setting $\widetilde{\mathcal{X}}((g,\gamma_h))=\mathcal{X}(g) \mathcal{X}(h)$ with $\gamma_h$ the inner automorphism associated to $h \in N$. We can easily see that $\widetilde{\mathcal{X}}$ and  $\widetilde{\mathcal{X}}^a$ are similar for all $a \in ((N \rtimes \Gamma_{\psi^\mathcal{H}}) \times \mathcal{H})_{\psi}$.
If $\mathcal{P}$ is now a representation of $G \rtimes \Gamma_\psi$ extending $\widetilde{\mathcal{X}}$, then the corresponding $\mu_a$ are characters of $(G \rtimes \Gamma_{\psi})/(G \rtimes \Inn(N))$ by Gallagher's lemma, see e.g. \cite[Corollary 6.17]{isaacs1994character}.
Further, note that the corresponding scalars on $C_{G\rtimes \Gamma_{\psi^\mathcal{H}}}(N)=\{(n,\gamma_n) \mid n \in N\}$ are all trivial.
\item It actually occurs that the character $\psi$ extends to $G \rtimes \Gamma_\psi$ but the extension cannot be chosen such that all $\mu_a$ are trivial.  An example for this is $\mathsf{G}_2(3)$ in its defining characteristic, see the proof of Proposition \ref{propnongeneric}.
\end{enumerate}
\end{rem}

\section{The groups $\mathsf{B}_2(2^i)$, $\mathsf{G}_2(3^i)$, $\mathsf{F}_4(2^i)$, and the Suzuki and Ree groups}\label{b2g2f4}
In this section, we verify the inductive McKay--Navarro condition for the groups with exceptional graph automorphisms and the Suzuki and Ree groups in their defining characteristic.
The groups $\mathsf{B}_2(2)'$, $\mathsf{G}_2(3)$, $\mathsf{F}_4(2)$, ${}^2\mathsf{B}_2(8)$, ${}^2\mathsf{G}_2(3)'$, and ${}^2\mathsf{F}_4(2)'$  will be studied separately in Section \ref{sectionnongeneric}, whereas ${}^2\mathsf{B}_2(2)$ is solvable and does not have to be considered. We follow \cite{maslowski2010} and \cite{ruhstorfer2017navarro} and extend the results from there.
\subsection{Notation} 
For integers $i \geq 2, f \geq 1$ ($f \geq 2$ for type $^2\mathsf{B}_2$), we let $$G \in \{\mathsf{B}_2(2^i), \mathsf{G}_2(3^i),\mathsf{F}_4(2^i),{}^2\mathsf{B}_2(2^{2f+1}), {}^2\mathsf{G}_2(3^{2f+1}), {}^2\mathsf{F}_4(2^{2f+1})\}$$ and let $p$ be the defining characteristic $2$ or $3$. We write $q:=p^{i}$ and $q^2:=p^{2f+1}$, respectively. The group $G$ is simple, non-abelian, has trivial Schur multiplier, and trivial centre \cite[Table 24.2, Remark 24.19]{malletesterman2011}. 

Let $\mathbf{G}$ be the corresponding simple algebraic group of type $\mathsf{B}_2,$  $\mathsf{G}_2$, or $\mathsf{F}_4$, respectively, defined over an algebraic closure $\mathbf{k}$ of $\mathbb{F}_p$. We denote by $F_p$ the standard Frobenius endomorphism of  $\mathbf{G}$ and by $\gamma$ the endomorphism of $\mathbf{G}$ induced by the exceptional graph automorphism $\rho$ of the Dynkin diagram.
If $G$ is a Suzuki or Ree group, we set $F=F_p^f \circ \gamma$;
otherwise $F=F_p^i$.
Then, $F$ is a Steinberg endomorphism such that $\mathbf{G}^F=G$. Note that the centre of $\mathbf{G}$ is also trivial by \cite[Theorem 1.12.5]{gorensteinlyonssolomon1998}. 

We fix an $F$-stable maximal torus $\mathbf{T} \subset \mathbf{G}$ and an $F$-stable Borel subgroup $\mathbf{B} \supset \mathbf{T}$. Let $\mathbf{U}$ be the unipotent radical of $\mathbf{B}$. By \cite[Corollary 24.11]{malletesterman2011}, $U=\mathbf{U}^F$ is a Sylow $p$-subgroup of $G$ with normalizer $B=\mathbf{B}^F$. Let $\Phi$ be the root system of $\mathbf{G}$ with respect to $\mathbf{T}$, $\Phi^\vee$ the set of coroots, $n$ the rank of $\Phi$, and $\Delta=\langle \alpha_i  \mid 1 \leq i \leq n \rangle$ the set of simple roots with respect to $\mathbf{B}$.

\subsection{Dual fundamental weights}\label{dualfundam}
 Let $X(\mathbf{T})$ be the character group of $\mathbf{T}$ and $Y(\mathbf{T})$ the group of cocharacters of $\mathbf{T}$. For $\alpha \in \Phi$, let $\mathbf{U}_\alpha$ be the root subgroup associated to $\alpha$ and fix an isomorphism $x_\alpha:(\mathbf{k},+) \rightarrow \mathbf{U}_\alpha$. 
Following \cite[Chapter 6]{maslowski2010}, we can now define dual fundamental weights $\omega_j^\vee \in Y(\mathbf{T})$ such that $\langle \alpha_i, \omega_j^\vee  \rangle = \delta_{ij}$ for $1 \leq i,j \leq n$.  
Since the centre of $\mathbf{G}$ is trivial, the dual fundamental weights are given by 
$$\omega_j^\vee: \mathbf{k}^\times \rightarrow \mathbf{T}, \quad \omega_j^\vee(\mu)= \prod_{k=1}^n h_{\alpha_k}(\mu^{la_{kj}})$$ 
with $(a_{ij})_{i,j} \in \mathbb{Q}^{n \times n}$ the inverse of the Cartan matrix of $\mathbf{G}$, $l$ the exponent of $X(\mathbf{T})/\mathbb{Z}\Phi$, and $h_{\alpha_k}(t)=x_{\alpha_k}(t)x_{-\alpha_k}(-t^{-1})x_{\alpha_k}(t)x_{\alpha_k}(-1)x_{-\alpha_k}(1)x_{\alpha_k}(-1)$ for all $t \in \mathbf{k}^\times$.
Note that they only depend on the type of the root system of $\mathbf{G}$ and not on $F$.

By definition, we have $\gamma^2=F_p$ and
$$\gamma(x_\alpha(t))=\left\{ \begin{array}{ll}
 x_{\rho(\alpha)}(t) & \text{if } \alpha \text{ is long}, \\
 x_{\rho(\alpha)}(t^p) & \text{if } \alpha \text{ is short,}
\end{array} \right.$$
see e.g. \cite[Theorem 1.15.4]{gorensteinlyonssolomon1998}. The outer automorphism group of $\mathbf{G}^F$ is cyclic with generator $\gamma$. If $\mathbf{G}^F$ is a Suzuki or Ree group, it has order $2f+1$ and is equal to $ \langle F_p \rangle$; otherwise it has order $2i$ \cite[Theorem 2.5.12]{gorensteinlyonssolomon1998}.  Since $\mathbf{U}=\prod_{\alpha \in \Phi^+} \mathbf{U}_\alpha$ with $\Phi^+$ the set of positive roots with respect to $\Delta$, $\gamma$ fixes $\mathbf{U}^F$. Thus, we have $\Gamma=\langle \gamma, \Inn(N) \rangle$.

We follow the ideas in \cite{ruhstorfer2017navarro} and generalize the considerations that were already made there for Frobenius endomorphisms and groups without exceptional graph automorphisms. 
\begin{lemma}\label{actiondualweights}
The action of $\gamma$ on the dual weights is given by
$$\gamma(\omega_j^\vee(\mu))=\left\{ \begin{array}{ll}
 \omega_{\rho(j)}^\vee(\mu) & \text{if } \alpha \text{ is long}, \\
 \omega_{\rho(j)}^\vee(\mu^p) & \text{if } \alpha \text{ is short,}
 \end{array} \right.$$ for $\mu \in \mathbf{k}^{\times}$ where we use $\rho$ to denote the permutation of the indices of the $\alpha_i$ induced by $\rho$.
\end{lemma}
\begin{proof} The claim follows from a simple computation using the inverses of the Cartan matrices as in \cite[Appendix]{maslowski2010}.
\end{proof}

\subsection{An equivariant bijection}
As already mentioned, 
we know that the inductive McKay condition holds in defining characteristic. More precisely, we know by \cite[Theorem 5]{brunat2009defining} that there exists a $\Gamma$-equivariant bijection 
$$\Omega_{\text{iMcK}}: {\Irr}_{p'}(\mathbf{G}^F) \rightarrow {\Irr}_{p'}(\mathbf{B}^F)$$ 
such that (2A) of Condition \ref{imck} holds.
We want to show that this bijection is also $\mathcal{H}$-equivariant.

Now, $(\mathbf{G},F)$ is selfdual by \cite[p.120]{carter1985finite} except for type $\mathsf{B}_2$ where we have a bijection between rational semisimple elements of $\mathbf{G}$ and its dual with an isomorphism of centralizers  \cite[p. 164]{lusztig1977}. In the following, we will use this bijection without further notice.
For a semisimple element $s \in {\mathbf{G}}^{F}$, we denote by $\mathcal{E}(\mathbf{G}^F, s)$ the set of irreducible constituents of all possible Deligne--Lusztig generalized characters $R_{\mathbf{T}'}^\mathbf{G}(s)$ for $F$-stable maximal tori $\mathbf{T}' \subseteq \mathbf{G}$. By the Jordan decomposition of irreducible characters of $\mathbf{G}^F$, we have a partition $$\Irr(\mathbf{G}^F)=\bigcup_{s}\mathcal{E}(\mathbf{G}^F, s)$$ where $s$ ranges over the conjugacy classes of semisimple elements of $\mathbf{G}^F$ \cite[Theorem 2.6.2]{geckmalle2020}.
We know how the outer automorphisms and $\mathcal{H}$ act on the character sets $\mathcal{E}(\mathbf{G}^F, s)$.

\begin{lemma} \label{lemmaHactionglobal}
Let $s \in \mathbf{G}^F$ be semisimple and $\psi \in \mathcal{E}(\mathbf{G}^F, s)$. For an outer automorphism $\kappa$ of $\mathbf{G}^F$, we have $\psi^\kappa \in \mathcal{E}(\mathbf{G}^F, s^\kappa)$.  
Further, for any $\sigma \in \mathcal{H}$ such that every $p'$-root of unity is mapped to its $p^k$-th power where $k$ is an integer, we have $\psi^\sigma \in \mathcal{E}(\mathbf{G}^F, s^{p^k})=\mathcal{E}(\mathbf{G}^F, s^{F_p^k}).$
\end{lemma}
\begin{proof}
The first statement is \cite[Proposition 1]{brunat2009defining} since all outer automorphisms of $\mathbf{G}$ commute with $F$.
Similar as \cite[Proposition 3.3.15]{geckmalle2020}, the second claim follows from the character formula and the fact that all maximal tori have an order prime to $p$.  
\end{proof}
We show an analogous result for the local case.
\begin{lemma} \label{lemmalocalcharHaction}
Let $\sigma \in \mathcal{H}$ such that every $p'$-root of unity is mapped to its $p^k$-th power where $k$ is an integer. Then, $\sigma$ acts on the elements of $\Irr_{p'}(\mathbf{B}^F)$ in the same way as $F_p^k$.
\end{lemma}
\begin{proof}
If $F$ is a Frobenius map, the claim is included in \cite[Remark 5.2]{ruhstorfer2017navarro}. We consider the Suzuki and Ree groups case by case. First, let $\mathbf{G}^F$ be of type ${}^2\mathsf{B}_2$ or ${}^2\mathsf{G}_2$. Then, $\mathbf{B}^F$ is a Frobenius group with $\mathbf{B}^F=\mathbf{U}^F \cdot \mathbf{T}^F$ where $\mathbf{U}^F \in {\Syl}_p(\mathbf{G}^F)$ and $\mathbf{T}^F$ is cyclic of order $q^2-1$ by \cite[Sect. XI.3]{huppertblackburn} and \cite[Proof of Lemma 5]{eaton2000dadesinductive}. Thus, the irreducible characters of $\mathbf{B}^F$ consist of the inflations of characters in $\Lin(\mathbf{T}^F)$ and the characters induced by $\Lin(\mathbf{U}^F) \setminus \{1_{\mathbf{U}^F}\}$. Since $F_p^k$ acts on the elements of $\mathbf{T}^ F$ by mapping them to their $p^k$-th power and $\mathbf{T}^ F$ has order prime to $p$, $F_p^k$ and $\sigma$ act in the same way on $\Lin(\mathbf{T}^F)$. This is inherited by the inflated characters. 

The non-trivial characters of $\mathbf{U}^F$ are all in the same orbit under $\mathbf{T}^F$ and induce the same unique irreducible character of $\mathbf{B}^F$ of degree $q^2-1$ \cite[(16C)]{isaacsmallenavarro2007} \cite[Proof of Lemma 5]{eaton2000dadesinductive}. This character has to be fixed under $F_p^k$ and $\sigma$ and it follows $\psi^{F_p^k}=\psi^{\sigma}=\psi$ for all $\psi \in \Irr(B)$.

Finally, let $\mathbf{G}^F$ be of type ${}^2\mathsf{F}_4$. With the notation and index sets as in \cite{himstedthuang2009charborel}, we have as in \cite[Proof of Lemma 6.1]{himstedthuang2009charborel}
$${\Irr}_{2'}(\mathbf{B}^F)=\{{}_B\chi_1(k,l)\}\cup \{{}_B\chi_2(k)\} \cup \{{}_B \chi_5(k)\} \cup \{{}_B\chi_8 \}.$$
As described there, $F_2$ acts on these characters by doubling the character parameters. By looking at the explicit character values given in \cite[Table A.6]{himstedthuang2009charborel}, we see that $\sigma$ acts on the characters in the same way as $F_2^k$.
\end{proof}

\begin{prop}\label{propdefiningequiv}
The $\Gamma$-equivariant bijection $\Omega_{\textrm{iMcK}}$ from \cite{brunat2009defining} is even $\mathcal{H}$-equivariant.
\end{prop}
\begin{proof}
Let $\sigma \in \mathcal{H}$. We know from \cite[Lemma 5]{brunat2009defining} that in every Lusztig series $\mathcal{E}(\mathbf{G}^F, s)$ there is only one character of $p'$-degree. Thus, $\sigma$ acts on $\Irr_{p'}(\mathbf{G}^F)$ and $\Irr_{p'}(\mathbf{B}^F)$ in the same way as some element $F_p^k \in \Gamma$. Since $\Omega_{\textrm{iMcK}}$ is $\Gamma$-equivariant, it follows for all $\psi \in \Irr_{p'}(G)$
$$\Omega_{\textrm{iMcK}}(\psi^\sigma)=\Omega_{\textrm{iMcK}}(\psi^{F_p^k })=\Omega_{\textrm{iMcK}}(\psi)^{F_p^k}=\Omega_{\textrm{iMcK}}(\psi)^{\sigma}.$$ This shows the claim.
\end{proof}

\subsection{Character extensions} 
We use the same construction and notation for the linear characters $\phi_S$ and the regular character $\xi$ of $\mathbf{U}^F/[\mathbf{U},\mathbf{U}]^F$ and $\mathbf{U}^F$ as in \cite[Section 3.2 and 4.4]{ruhstorfer2017navarro} and just briefly recall the properties we need later. Note that the construction still works for Steinberg maps $F$. 
We set $r=n$ if $F$ is a Frobenius map and $r=n/2$ if we are in the case of Suzuki and Ree groups.
For $1 \leq i \leq r$, we have an isomorphism 
$$x_i:(\mathbb{F}_{q^{n/r}},+) \rightarrow  \prod_{\alpha \in \rho(\{\alpha_i\})} \mathbf{U}_{\alpha}$$ where $\rho$ is the automorphism on the roots induced by $F$ as before. We choose a non-trivial character $\phi_0 \in \Irr(\mathbb{F}_{q^{n/r}},+)$ that is invariant under field automorphisms which is possible by \cite[Lemma 11.12]{maslowski2010}. For $S \subseteq \{1, \ldots , r \}$, the character
$\phi_S \in \Irr(\mathbf{U}/[\mathbf{U}, \mathbf{U}]^F) $ is given by
$$\phi_S(x_1(a_1) \ldots x_r(a_r))= \prod_{i \in S} \phi_0(a_i)$$ and we write $\xi:=\phi_{\{1, \ldots , r\}}$. By construction, the characters $\phi_S$ are invariant under field automorphisms. If $\rho$ is trivial, then the isomorphisms $x_i$ are just the root homomorphisms $x_\alpha$ and we can easily determine the action of $\gamma$ on the $\phi_S$. We see that $\xi$ is always fixed by $\langle \gamma \rangle$.

The character $\Gamma_1=\Ind_{\mathbf{U}^F}^{\mathbf{G}^F} (\xi)$ is the Gelfand--Graev character of $\mathbf{G}^F$, see \cite[Paragraph before Definition 12.3.3]{dignemichel2020representations}. 
We denote the Alvis--Curtis duality map for $\mathbf{G}$ by $D_{\mathbf{G}}: \mathbb{Z}\Irr(\mathbf{G}^F) \rightarrow \mathbb{Z}\Irr(\mathbf{G}^F) $, see \cite[Section 7.2]{dignemichel2020representations} or \cite[Section 3.4]{geckmalle2020}.
\begin{lemma}\label{lemmat}
For every $\sigma \in \mathcal{G}$ there exists $t \in \mathbf{T}^{ \gamma}$ such that $\phi_S^\sigma=\phi_S^t$ for all $S\subseteq \{1,\ldots , r\}$.
\end{lemma}
\begin{proof}
Let $\phi_0 \in \Irr(\mathbb{F}_{q^{n/r}},+)$ be as in \cite[Section 3.2]{ruhstorfer2017navarro}. As in the proof of \cite[Lemma 6.3(a)]{ruhstorfer2017navarro}, let $b \in \mathbb{F}_p^\times$ be such that $\phi_0^\sigma(a)=\phi_0(ba)$ for all $a \in \mathbb{F}_{q^{n/r}}$. It is shown there that for $s_i:=\omega_i^\vee(b) \in \mathbf{T}^{F_p}$ with $1 \leq i\leq n$ and $t:=\prod_{i=1}^n s_i \in \mathbf{T}^{F_p}$ we have $\phi_S^t=\phi_S^\sigma$. Since $b^p=b$, it follows with Lemma \ref{actiondualweights} that $$\gamma(t)=\prod_{i=1}^n \gamma(s_i)=\prod_{i=1,\alpha_i \textrm{ long}}^n \omega_{\rho(i)}^\vee(b) \prod_{i=1,\alpha_i \textrm{ short}}^n \omega_{\rho(i)}^\vee(b^p)
=\prod_{i=1}^n \omega_{i}^\vee(b)=t.$$ Thus, $t$ is $\gamma$-invariant.
\end{proof}
We now want to show that the global and local characters can be extended to their stabilizers in $\Gamma$ such that the corresponding $\mu_a$ are trivial.
We generalize the constructions from \cite[Section 6.3]{ruhstorfer2017navarro} to our case, see also \cite{digne1999descenteshintani} for more details. Let $F_0$ be a Steinberg automorphism of $\mathbf{G}$ such that $F_0^k=F$ for some integer $k$. For any $F$-stable closed subgroup $\mathbf{H}$ of $\mathbf{G}$ we set
$$\underline{\mathbf{H}}:=\mathbf{H} \times F_0^{k-1}(\mathbf{H}) \times \ldots \times F_0(\mathbf{H})\leq \mathbf{G}^k=:\underline{\mathbf{G}}.$$
We define the automorphism $$\tau: \underline{\mathbf{G}} \rightarrow \underline{\mathbf{G}}, \quad \tau(g_1,\ldots, g_k)=(g_2, \ldots , g_k,g_1).$$
As in \cite[Section 6.3]{ruhstorfer2017navarro}, we see that $\tau$ is quasi-central as defined in \cite[D\'efinition-Th\'eor\`eme 1.15]{dignemichel1994grpsnonconn}. The projection onto the first coordinate $\pr$ induces isomorphisms 
$$\underline{\mathbf{H}}^{\tau F_0}\cong \mathbf{H}^F , \quad \underline{\mathbf{G}}^{\tau F_0} \rtimes \langle \tau \rangle \cong \mathbf{G}^{F} \rtimes \langle F_0 \rangle  .$$
The following proposition extends \cite[Proposition 6.7]{ruhstorfer2017navarro} to our setting. Although Ruhstorfer considers only Galois automorphisms in $\mathcal{H}$, the proof given there also holds for all Galois automorphisms that satisfy \cite[Assumption 6.6]{ruhstorfer2017navarro} and act like some field automorphism of the group. Thus, the changes we make are only due to the generalization to Steinberg automorphisms $F_0$ and the limited knowledge about duality for disconnected groups in this setting.
\begin{prop} \label{propsemisimpleextension}
Let $\psi \in \Irr(\mathbf{G}^F)$ be a semisimple character and $\sigma \in \mathcal{G}$ a Galois automorphism such that we have $\psi^\sigma=\psi^{F_p^e}$ for some $e \in \mathbb{N}$. 
Then, there exists an extension $\hat{\psi} \in \Irr(\mathbf{G}^F \rtimes \langle \gamma \rangle_\psi)$ of $\psi$ such that $\hat{\psi}^{F_p^e t \sigma^{-1}}=\hat{\psi}$ where $t \in \mathbf{T}^\gamma$ is as in Lemma \ref{lemmat}.
\end{prop}
\begin{proof}
Since $\langle \gamma \rangle_{\psi}$ is cyclic and $\psi^F=\psi$, we find a generator $F_0$ of $\langle \gamma \rangle_{\psi}$ such that $F_0^k=F$ for some $k \in \mathbb{N}$. 

\textbf{Step 1: Translating the characters to $\mathbf{\underline{G}}$.}
Denote by $${\pr}_{\mathbf{G}\rtimes \langle \tau \rangle }^{\vee}: \mathbb{Z}\Irr(\mathbf{G}^{F} \rtimes \langle F_0 \rangle) \rightarrow \mathbb{Z}\Irr( \underline{\mathbf{G}}^{\tau F_0} \rtimes \langle \tau \rangle ),$$ $${\pr}_{\mathbf{G}}^{\vee}: \mathbb{Z}\Irr(\mathbf{G}^F) \rightarrow \mathbb{Z}\Irr(\underline{\mathbf{G}}^{\tau F_0}), \quad {\pr}_{\mathbf{U}}^{\vee}: \mathbb{Z}\Irr(\mathbf{U}^F) \rightarrow \mathbb{Z}\Irr(\underline{\mathbf{U}}^{\tau F_0})$$ the maps induced by $\pr$, restricting to isomorphisms of the respective sets of irreducible characters. 
Note that the centers of $\mathbf{G}$ and $\underline{\mathbf{G}}$ are trivial. We set $\underline{\psi}:={\pr}_{\mathbf{G} }^\vee(\psi)$, $ \underline{\xi}:=\pr_{\mathbf{U}}^{\vee}(\xi)$, and $\underline{\Gamma_1}=\Ind_{\underline{\mathbf{U}}^{\tau F_0}}^{\underline{\mathbf{G}}^{\tau F_0}}(\underline{\xi})$. Using the character formula for induction, we see that $${\pr}_{\mathbf{G}}^\vee \circ {\Ind}_{\mathbf{U}^F}^{\mathbf{G}^F}={\Ind}_{\underline{\mathbf{U}}^{\tau F_0}}^{\underline{\mathbf{G}}^{\tau F_0}} \circ {\pr}_{\mathbf{U}}^{\vee}$$ and it follows $\underline{\Gamma_1}={\pr}_{\mathbf{G}}^\vee(\Gamma_1)$. 

Since $\psi$ is semisimple, it is a constituent of $D_{{\mathbf{G}}}(\Gamma_1)$ by \cite[Corollary 12.4.10]{dignemichel2020representations}. From \cite[Proposition 3.4.3]{geckmalle2020} we know that Alvis--Curtis duality and $\pr_{\mathbf{G}}^\vee$ also commute in this setting. Consequently, we have $D_{\underline{\mathbf{G}}} \circ \pr_{\mathbf{G}}^\vee= {\pr}_{\mathbf{G}}^\vee \circ D_{\mathbf{G}}$ and $\underline{\psi}$ is a constituent of $D_{\underline{\mathbf{G}}}(\underline{\Gamma_1})$.
Let $\phi:\mathbf{G} \rightarrow \mathbf{G}$ be the morphism given by the action of $F_p^e t$ and set $$\underline{\phi}: \underline{\mathbf{G}} \rightarrow \underline{\mathbf{G}}, \quad (g_1, \ldots g_k) \mapsto (\phi(g_1),\ldots \phi(g_k)).$$ Then, we have $\underline{\psi}^{\underline{\phi}\sigma^{-1}}=\underline{\psi}$. By the choice of $t$ and since $\xi$ is invariant under field automorphisms, it follows that $\xi^{\phi\sigma^{-1}}=\xi^{t\sigma^{-1}}=\xi$. This implies $\underline{\xi}^{\underline{\phi}\sigma^{-1}}=\underline{\xi}$ and $\underline{\Gamma_1}^{\underline{\phi}\sigma^{-1}}=\underline{\Gamma_1}$.

\textbf{Step: 2 Extending characters to the disconnected group.}
Since $\xi$ is fixed by $F_0$, $\underline{\xi}\in \Irr( \underline{\mathbf{U}}^{\tau F_0})$ is $\tau$-invariant and we can extend $\underline{\xi} $ to a character  $\hat{\underline{\xi}} \in \Irr(\underline{\mathbf{U}}^{\tau F_0} \rtimes \langle \tau \rangle)$ by setting $\hat{\underline{\xi}}(\tau)=1$. Let $$\underline{\hat{\Gamma}_1}:= {\Ind}_{\underline{\mathbf{U}}^{\tau F_0}\rtimes \langle \tau \rangle}^{\underline{\mathbf{G}}^{\tau F_0}\rtimes \langle \tau \rangle}(\hat{\underline{\xi}}) .$$ Then, $\underline{\hat{\Gamma}_1}$ extends $\underline{\Gamma_1}$ by the Mackey formula for characters \cite[Theorem 1.16]{navarro2018book}. 

We need the construction of Deligne--Lusztig characters and duality for disconnected groups from \cite[Definition 2.2 and 3.10]{dignemichel1994grpsnonconn}. There, Digne and Michel only consider finite groups of Lie type that arise via Frobenius endomorphisms. Since $\tau F_0$ is not necessarily a Frobenius endomorphism of $\underline{\mathbf{G}} \rtimes \langle \tau \rangle$, we have to make sure that the constructions still work in our setting. Going through \cite{dignemichel1994grpsnonconn}, it is easy to see that most of the claims still hold if we allow Steinberg endomorphisms. We checked this for all results and definitions that we need in the following.
Using \cite[Corollaire 2.4.]{dignemichel1994grpsnonconn}, it is easy to see that 
$${\Res}_{\underline{\mathbf{G}}^{\tau F_0}}^{(\underline{\mathbf{G}} \rtimes \langle \tau \rangle )^{\tau F_0}} \circ D_{\underline{\mathbf{G} } \rtimes \langle \tau \rangle }=D_{\underline{\mathbf{G}}} \circ  {\Res}_{\underline{\mathbf{G}}^{\tau F_0}}^{(\underline{\mathbf{G}} \rtimes \langle \tau \rangle )^{\tau F_0}}.$$
Thus, $D_{\underline{\mathbf{G}} \rtimes \langle \tau \rangle }(\underline{\hat{\Gamma}_1})$ extends $D_{\underline{\mathbf{G}}}(\underline{\Gamma_1})$.

\textbf{Step 3: Conclusion.}
Here, Harish-Chandra induction and restriction still commute with Galois automorphisms and group automorphisms in the generalized case. Thus, we have
$$D_{\underline{\mathbf{G}} \rtimes \langle \tau \rangle}(\underline{\hat{\Gamma}_1})^{\underline{\phi} \sigma^{-1}}=D_{\underline{\mathbf{G}} \rtimes \langle \tau \rangle}(\underline{\hat{\Gamma}_1}^{\underline{\phi} \sigma^{-1}})=D_{\underline{\mathbf{G}} \rtimes \langle \tau \rangle}(\underline{\hat{\Gamma}_1}).$$
As mentioned before, the character $\underline{\psi}$ is a constituent of $D_{\underline{\mathbf{G}}}(\underline{\Gamma_1})$. Since $D_{\underline{\mathbf{G}}}(\underline{\Gamma_1})$ and $D_{\underline{\mathbf{G}} \rtimes \langle \tau \rangle}(\underline{\hat{\Gamma}_1})$ are both multiplicity free, there is only one constituent $\hat{\underline{\psi}}$ of $D_{\underline{\mathbf{G}} \rtimes \langle \tau \rangle}(\underline{\hat{\Gamma}_1})$ that extends $\underline{\psi}$. Thus, it is $\underline{\phi}\sigma^{-1}$-stable and it follows that $({\pr_{\mathbf{G} \rtimes \langle \tau \rangle}^{\vee}})^{-1}(\hat{\underline{\psi}}) \in  \Irr(\mathbf{G}^F \rtimes \langle F_0 \rangle)$ is a ${\phi}\sigma^{-1}$-stable character extending $\psi$.
\end{proof}
It remains to show the existence of suitable extensions of the local characters.
\begin{prop} \label{proplocaldefiningextension}
Let $\chi \in \Irr(\mathbf{B}^F)$ be a semisimple character and $\sigma \in \mathcal{G}$ a Galois automorphism such that we have $\chi^\sigma=\chi^{F_p^e}$ for some $e \in \mathbb{N}$.  Let $t \in \mathbf{T}^\gamma$ as is Lemma \ref{lemmat} and set $x_\sigma:=F_p^et\sigma^{-1}$. Then, there exists an extension $\hat{\chi} \in \Irr(\mathbf{B}^F \rtimes \langle \gamma \rangle_\psi)$ of $\chi$ such that $\hat{\chi}^{x_\sigma}=\hat{\chi}$. 
\end{prop}
\begin{proof}
 For the groups $\mathsf{B}_2(2^i), \mathsf{G}_2(3^i)$, and $\mathsf{F}_4(2^i)$, this can be shown in the same way as in \cite[Prop. 6.10]{ruhstorfer2017navarro}. For the Suzuki and Ree groups, we give a different proof. Let $\mathbf{G}^F$ be ${}^2\mathsf{B}_2(q^2)$ or ${}^2\mathsf{G}_2(q^2)$. We recall the $p'$-characters of $\mathbf{B}^F$ from the proof of Lemma \ref{lemmalocalcharHaction}:
$${\Irr}_{p'}(\mathbf{B}^F)=\{{\Inf}_{\mathbf{T}^F}^{\mathbf{B}^F} (\tau_i) \mid \tau_i \in \Lin(\mathbf{T}^F)\} \cup \{{\Ind}_{\mathbf{U}^F}^{\mathbf{B}^F}(\xi) \}.$$
Since the characters $\chi_i:=\Inf_{\mathbf{T}^F}^{\mathbf{B}^F} (\tau_i)$ are linear, we obtain an extension to $\mathbf{B}^F \rtimes \langle \gamma \rangle_{\chi_i}$ that is invariant under $F_p^e\sigma^{-1}$ by setting $\hat{\chi_i}(b \gamma^k)=\chi_i(b)$ for all $k \in \mathbb{N}$, $b \in \mathbf{B}^F $. These extensions are also invariant under $x_\sigma$ since $t \in \mathbf{B}^F $.

We now consider the linear character $\xi \in \Irr(\mathbf{U}^F )$. As noted before, $\xi$ is invariant under $x_\sigma$ and can be extended to $\hat{\xi} \in \Irr(\mathbf{U}^F \rtimes  \langle \gamma \rangle)$ such that it is $x_\sigma$-invariant. This also gives us an $x_\sigma$-invariant extension ${\Ind}_{\mathbf{U}^F \rtimes \langle \gamma \rangle}^{\mathbf{B}^F \rtimes \langle \gamma \rangle}(\hat{\xi}) $ of the character ${\Ind}_{\mathbf{U}^F}^{\mathbf{B}^F}(\xi) $.

For $\mathbf{G}^F={}^2\mathsf{F}_4(q^2)$, the characters of $\mathbf{B}^F$ were explicitly constructed in \cite{himstedthuang2009charborel}. With the notation and suitable index sets given there, we have 
$$\Irr_{2'}(\mathbf{B}^F)=\{ {}_B\chi_1(k,l)\} \cup \{{}_B\chi_2(k) \} \cup \{{}_B\chi_5(k) \} \cup \{{}_B\chi_8 \}.$$
Since ${}_B\chi_1(k,l)$ is linear, we can extend it as claimed. We write $H:=C_{\mathbf{T}^F}(\alpha_1(1))\mathbf{U}^F$ and
denote by $\lambda_2(k) \in \Irr(H)$ the character inducing ${}_B\chi_2(k)$ as given in \cite[p. 9]{himstedthuang2009charborel}. Then, we can easily see $\lambda_2(k)^{F_2}=\lambda_2(2k)$ and it follows 
$$D:=\langle F_2 \rangle_{{}_B\chi_2(k)}=\langle F_2 \rangle_{\lambda_2(k)}.$$ The character $\lambda_2(k)$ restricts to a linear character $\phi_{\{1\}} \in \Irr(\mathbf{U}^F)$. The action of $\mathbf{T}^F$ on the root subgroups is described in \cite[Table 2]{himstedthuang2009charborel} and $H$ is the inertia subgroup of $\phi_{\{1\}}$ in $\mathbf{B}^F$. By Clifford correspondence, $\lambda_2(k)$ is the unique character of $H$ inducing ${}_B\chi_2(k)$ such that $\lambda_2(k)|_{\mathbf{U}^F}$ has $\phi_{\{1\}}$ as a constituent. 
Since $\phi_{\{1\}}$ and ${}_B\chi_2(k)$ are $x_\sigma$-invariant, the Clifford correspondent $\lambda_2(k)$ is also $x_\sigma$-invariant. Thus, we can extend the linear character $\lambda_2(k)$ to an $x_\sigma$-invariant character $\hat{\lambda}_2(k) \in \Irr(H \rtimes D)$. 
It follows that $\Ind^{\mathbf{B}^F \rtimes D}_{H \rtimes D}(\hat{\lambda}_2(k))$ is an $x_\sigma$-invariant extension of ${}_B\chi_2(k)$.
In the same way, this can be shown for the other characters. 
\end{proof}

\subsection{Verification of the inductive condition} We are now able to verify the inductive McKay--Navarro condition for all considered groups in their defining characteristic.
\begin{theorem}
The groups $\mathsf{B}_2(2^i), \mathsf{G}_2(3^i),
\mathsf{F}_4(2^i), {}^2\mathsf{B}_2(2^{2f+1}), {}^2\mathsf{G}_2(3^{2f+1}), $ ${}^2\mathsf{F}_4(2^{2f+1})$ satisfy Condition \ref{imck} with $p$ the defining characteristic $p=2$ or $p=3$, respectively.
\end{theorem}
\begin{proof}
We use the notation of Condition \ref{imck} and let $S$ be one of the groups above. As already mentioned, the group $S$ is simple, non-abelian and has trivial Schur multiplier; thus we can consider $G=S$. We want to verify the condition for $N=B$.

Let $\Omega$ be the $\Gamma \times \mathcal{H}$-equivariant bijection from Proposition \ref{propdefiningequiv} and $D=\langle \gamma \rangle$. 
We know by Lemma \ref{lemmaHactionglobal} that every $\sigma \in \mathcal{H}$ acts on $\psi \in \Irr_{p'}(G)$ in the same way as $F_p^e$ for some $e \in \mathbb{N}$. Therefore, every  $a \in ( D \times \mathcal{H})_\psi$ is of the form $(dF_p^{e}, \sigma^{-1})$ for some $d \in D_\psi$. As before, let $x_\sigma:=F_p^et\sigma^{-1}$ with $t \in \mathbf{T}^\gamma$ as in Lemma \ref{lemmat}. 
We know by Proposition \ref{propsemisimpleextension} that we find an extension $\hat{\psi} \in \Irr(G \rtimes D_\psi )$ such that $\hat{\psi}^{x_\sigma}=\hat{\psi}$ for all $\sigma \in \mathcal{H}$. 
As we have  $t \in \mathbf{T}^D  \subseteq B \subseteq G$ and $d \in D_\psi$, both act trivially on $\hat{\psi}$ and we have $$\hat{\psi}=\hat{\psi}^{x_\sigma}=\hat{\psi}^{F_p^e \sigma^{-1}}=\hat{\psi}^{dF_p^e \sigma^{-1}}.$$
We can now extend $\hat{\psi}$ to the inner automorphisms in $\Gamma$ as described in Remark \ref{rem1}(\ref{remgallagher}). If $\mathcal{P}$ is a representation affording this extended character, then $\mathcal{P}$ and $\mathcal{P}^a$ afford the same character for all $a \in ((B \rtimes \Gamma)\times \mathcal{H})_\psi$ and are thereby similar. Therefore, the $\mu_a$ from Lemma \ref{lemmamu_a}(\ref{defmu_a}) are trivial. Using Proposition \ref{proplocaldefiningextension}, the same can be done for the local character $\Omega(\psi)$ and the claim follows.
\end{proof}

\section{The groups $\mathsf{B}_n(2)$}\label{bn2}
We verify Condition \ref{imck} for $S=\mathsf{B}_n(2)$ with an integer $n \geq 4$ and $p=2$. For $n=2$ and $n=3$, this will be treated separately in Proposition \ref{propnongeneric}.
\subsection{Action of Galois automorphisms in the global case}
The group $S$ is simple, non-abelian and has trivial Schur multiplier \cite[Remark 24.19]{malletesterman2011}, hence we can consider Condition \ref{imck} for $G=S$. Note that the connected reductive group $\mathsf{B}_n$ defined over $\overline{\mathbb{F}_2}$ has connected centre. 
\begin{lemma}\label{lemmaHG}
The Galois automorphisms in $\mathcal{H}$ act trivially on ${\Irr}_{2'}(G)$.
\end{lemma}
\begin{proof}
By \cite[p. 164]{lusztig1977}, there exists a bijection between the rational semisimple elements of ${\Sp}_{2n}(2)$ and of its dual $\SO_{2n+1}(2)$ with an isomorphism of centralizers. Thus, the Jordan decomposition of $\chi \in \Irr_{2'}(G)$ can be written as $(s, \nu)$ with $s \in G$ semisimple and $\nu$ a unipotent character of $C_{G}(s)$. 

Let $\sigma \in \mathcal{H}$ be such that every $2'$-th root of unity is mapped to its $2^k$-th power for some $k \in \mathbb{Z}$. Since the semisimple elements of $G$ have odd order, $\chi^\sigma$ has Jordan decomposition $(s^{2^k}, \nu^\sigma)$ by \cite{srinivasanvinroot2019}. As in the proof of \cite[Proposition 2]{cabanesoddchardeg}, we have $C_G(s) \cong \Sp_{2j}(2) \times C$ for some $0 \leq j \leq n$ and $C$ a product of finite linear or unitary groups. By \cite[Corollary 1.12]{lusztig2002rat}, every unipotent character of a group of type $\mathsf{A}$ or $\mathsf{B}$ is rational-valued and it follows that $\nu^\sigma=\nu$. 

It remains to show that $s$ and $s^{2^k}$ are conjugate for all $k \in \mathbb{Z}$.
The conjugacy classes of $G$ are uniquely determined by the characteristic polynomials of their elements. Let $\pi=X^n+a_{n-1}X^{n-1} + \cdots + a_0 \in \mathbb{F}_2[X]$ be the characteristic polynomial of $s$ with roots $ \lambda_1, \ldots \lambda_n$ over $\overline{\mathbb{F}_2}$. By Vieta's formula, we have 
$a_{n-i}=e_i(\lambda_1, \ldots , \lambda_n)$ where $e_i$ is the $i$-th elementary symmetric polynomial in $n$ variables. 
Applying the Frobenius map on $\overline{\mathbb{F}_2}$ yields $a_{n-i}=e_i(\lambda_1^2, \ldots , \lambda_n^2)$ since all coefficients of the $e_i$ lie in $\mathbb{F}_2$. 

We know that the eigenvalues of $s^2$ over $\overline{\mathbb{F}_2}$ are given by $ \lambda_1^2, \ldots \lambda_n^2$. Again by Vieta's formula, the coefficients of the characteristic polynomial of $s^2$ coincide with the coefficients of $\pi$. Therefore, $s$ and $s^{2}$ are conjugate wich implies that $s$ and $s^{2^k}$ are conjugate for all $k \in \mathbb{Z}$. This shows that $\chi^\sigma=\chi$.
\end{proof}
\subsection{Action of Galois automorphisms in the local case}
Let $\SymM_n(2)$ be the additive group of symmetric $n \times n$ -matrices over $\mathbb{F}_2$, ${U}_n(2) \leq \GL_n(2)$ the group of upper triangular unipotent matrices over $\mathbb{F}_2$ and $R$ a Sylow $2$-subgroup of $G$.
By \cite[Proposition 3]{cabanesoddchardeg}, $R$ is self-normalizing. We want to show that Condition \ref{imck} is satisfied for $N=R$.
\begin{lemma}\label{lemmaHR}
The Galois automorphisms in $\mathcal{H}$ act trivially on ${\Irr}_{2'}(R)$.
\end{lemma}
\begin{proof}
As in \cite{cabanesoddchardeg}, we see $R \cong {\SymM}_n(2) \rtimes {U}_n(2)$ with ${U}_n(2)$ acting on ${\SymM}_n(2)$ by $x.s=xsx^T$ for $x \in {U}_n(2)$ and $s \in {\SymM}_n(2)$. With the same considerations as in \cite[proof of Proposition 3]{cabanesoddchardeg}, it follows
\begin{align*}
R/R' &\cong 
({\SymM}_n(2)/[{\SymM}_n(2),{U}_n(2)]) \times ({U}_n(2)/{U}_n(2)') \cong (C_{2})^{n+1}.
\end{align*}
Thus, every value of a linear character of $R/R'$ is either $1$ or $-1$ and thereby an integer. The claim follows since all linear characters of $R$ can be obtained as inflated characters of $R/R'$.
\end{proof}
\begin{prop}
Condition \ref{imck} is satisfied for the group $S=\mathsf{B}_n(2)$ with $n\geq 4$ and $p=2$.
\end{prop}
\begin{proof}
The group $G=S$ has trivial outer automorphism group \cite[Section 2.5]{gorensteinlyonssolomon1998}. Thus, $\Gamma= \Inn(R)$ acts trivially on all characters in ${\Irr}_{2'}(G)$ and ${\Irr}_{2'}(R)$. 
By Lemma \ref{lemmaHG} and Lemma \ref{lemmaHR}, the group $\mathcal{H}$ also acts trivially on these characters. We know  $|{\Irr}_{2'}(G)|=|{\Irr}_{2'}(R)|=2^{n+1}$ by \cite[Section 2]{cabanesoddchardeg}. Thus, there obviously exists a $\Gamma \times \mathcal{H}$-equivariant bijection between ${\Irr}_{2'}(G)$ and ${\Irr}_{2'}(R)$. We can extend all characters to the inner automorphisms as in Remark \ref{rem1} (\ref{remgallagher}) and see that the extension part of Condition \ref{imck} is also satisfied. This shows the claim.
\end{proof}

\section{Groups with non-generic Schur multiplier, $\mathsf{G}_2(2)'$, $^2\mathsf{G}_2(3)'$, and $^2\mathsf{F}_4(2)'$} \label{sectionnongeneric}
In this section, we show that the inductive McKay--Navarro condition is also satisfied for $\mathsf{G}_2(2)'$, $^2\mathsf{G}_2(3)'$, $^2\mathsf{F}_4(2)'$, and the finite groups of Lie type with non-generic Schur multiplier (see e.g. \cite[Table 24.3]{malletesterman2011}) in their defining characteristic.

\subsection{About $p$-extensions} 
As we see in the following lemma, we often do not have to consider the full Schur cover. Given a group $S$ and some covering group $G$ of $S$, the $p'$-part of the covering group is some covering group $H$ with $G \twoheadrightarrow H \twoheadrightarrow S$ such that $|G|/|H|$ is a $p$-power and $|H|/|S|$ is prime to $p$.
\begin{lemma}\label{lempext}
Let $S$ be a simple non-abelian group, $p$ a prime, $G$ the Schur cover of $S$ and $H$ the $p'$-part of the Schur cover. Let $M$ be the normalizer of a Sylow $p$-subgroup $P$ of $H$ and set $\Gamma=\Aut(H)_P$. Assume that there exists a $\Gamma \times \mathcal{H}$-equivariant bijection 
$\Omega': {\Irr}_{p'}(H) \rightarrow {\Irr}_{p'}(M)$ such that for all $\psi \in \Irr_{p'}(H)$ we have
$$(H \rtimes \Gamma_{\psi^\mathcal{H}}, H, \psi) \geq_c (M \rtimes \Gamma_{\psi^\mathcal{H}}, M, \Omega'(\psi)).$$ Then the inductive McKay--Navarro condition holds for $S$.
\end{lemma}
\begin{proof}
We have $G/Z(G) \cong S$ and $G/Z(G)_{p} \cong H$; thus $G$ is a central $p$-extension of $H$. 
Consider an irreducible character of $G$ with degree $d$ afforded by a representation $\mathcal{X}$ such that its restriction to $Z(G)_{p}$ is not trivial. Then we find a $z \in Z(G)_{p}$ such that $\mathcal{X}(z)=\varepsilon I_d$ for some $\varepsilon \in \mathbb{C}^\times$ with $\varepsilon^{|Z(G)_{p}|}=1$ and $\varepsilon \neq 1$. Since $S$ is perfect, $G$ is also perfect and its only linear character is the trivial one. Thus, the determinant of $\mathcal{X}(g)$ has to be $1$ for all elements $g \in G$ and it follows  $$\det (\mathcal{X}(z))=\varepsilon^d=1.$$ Since $|Z(G)_{p}|$ is a $p$-power and $\varepsilon \neq 1$, $d$ and $p$ cannot be coprime and it follows $p \mid d$. 
Therefore, all irreducible $p'$-characters of $G$ have $Z(G)_{p}$  in their kernel and can be obtained by inflating irreducible $p'$-characters of $H$. 
Thus, deflation yields a bijection $\Def^G_H$ between $p'$-characters of $G$ and $H$.

Let $R \leq G$ be the preimage of $P$ and $N$ its normalizer in $G$. 
Since $Z(G)_{p} \leq N$ by definition and Sylow theory, $N$ is a central $p$-extension of $M$. 
Thus, inflation yields a bijection $\Inf^N_M$ between the $p'$-characters of $M$ and $N$. It follows that $\Omega'$ induces a bijection $$\Omega:= {\Inf}^N_M \circ \Omega' \circ {\Def}^G_H: {\Irr}_{p'}(G) \rightarrow {\Irr}_{p'}(N).$$
By \cite[Corollary B.8]{navarro2018book}, the automorphism groups of $G$ and $H$ are both isomorphic to the automorphism group of $S$. Consequently, every automorphism of $H$ stabilizing $P$ extends uniquely to an automorphism of $G$ stabilizing $R$ and we can consider $\Gamma$ as a subset of $\Aut(G)$. 
Since every $\chi \in {\Irr}_{p'}(G)$ is uniquely determined by the values of $\chi$ on $H$, $\Gamma$ and $\mathcal{H}$ act on ${\Irr}_{p'}(G)$ in the same way as on ${\Irr}_{p'}(H)$, i.e. ${\Def}^G_H$ is $\Gamma \times \mathcal{H}$-equivariant. The same is true for $\Inf^N_M$; thus $\Omega$ is $\Gamma \times \mathcal{H}$-equivariant. 

It remains to show that the extension conditions in \ref{imck} are satisfied.
Given projective representations of $H \rtimes \Gamma_{\Def^G_H(\psi)}$ and $M \rtimes \Gamma_{\Omega'(\Def^G_H(\psi))}$ such that (2A) and (2B) of  Condition \ref{imck} hold, we can extend them trivially to $Z(G)_{p}$. These extended projective representations still satisfy (2A) and (2B) and Condition \ref{imck} holds for $S$.
\end{proof}
\subsection{Explicit computations}
We first state an easy consequence of \cite[(6.28)]{isaacs1994character} that sometimes ensures the existence of suitable character extensions.
\begin{lemma} \label{lemmaextensiontrivial}
In the setting of Condition \ref{imck}, let $\Gamma'$ be a subgroup of $\Gamma$ containing representatives for all outer automorphisms in $\Gamma$. Assume that $\Gamma'$ is a $p$-group. Then, for every $\psi \in \Irr_{p'}(G)$ there exists a $(\Gamma \times \mathcal{H})_\psi$-invariant extension of $\psi$ to $G \rtimes \Gamma_\psi$ such that all scalars associated to $C_{G \rtimes \Gamma}(N)$ are $1$.
\end{lemma}
\begin{proof}
The determinant of $\psi$ is trivial because $G$ is perfect. Since $\Gamma'$ is a $p$-group, it is solvable and $(|\Gamma'|,\psi(1))=1$. By \cite[(6.28)]{isaacs1994character}, we find a unique extension $\hat{\psi} \in \Irr(G \rtimes \Gamma'_\psi)$ such that $\hat{\psi}$ has trivial determinant. We easily see $\det(\hat{\psi}^{(x,\sigma)})=\det(\hat{\psi})$ for all $(x, \sigma) \in (\Gamma \times \mathcal{H})_\psi$ and with the uniqueness we can conclude $\hat{\psi}^{(x,\sigma)}=\hat{\psi}$. Now let $g \in C_{G \rtimes \Gamma}(N)$ and $\varepsilon \in \mathbb{Q}^{\ab}$ be the associated scalar, i.e. $\hat{\psi}(g)=\varepsilon \psi(1)$. Then, we have $\det(\hat{\psi}(g))=\varepsilon^{\psi(1)}=1$ and $\varepsilon^{\ord(g)}=1$. Since $\ord(g)$ and $\psi(1)$ are coprime, it follows $\varepsilon=1$ and we get the claim by extending the characters canonically to the inner automorphisms in $\Gamma$.
\end{proof}

We now verify the inductive McKay--Navarro condition for the finitely many groups that were not considered so far.
In the following proof, all explicit computations were made with \textsf{GAP} \cite{GAP4}.
\begin{prop}\label{propnongeneric}
The inductive McKay--Navarro condition holds for $\mathsf{G}_2(2)'$, $^2\mathsf{G}_2(3)'$, $^2\mathsf{F}_4(2)'$, and the simple groups of Lie type with non-generic Schur multiplier in their defining characteristic.
\end{prop}
\begin{proof}
Let $S$ be a simple group of Lie type, $p$ its defining characteristic and $G$ the $p'$-part of the Schur cover of $S$. Note that the proof of \cite[Theorem 7.3]{ruhstorfer2017navarro} does not use the fact that $\mathbf{G}^F=G$ is the universal covering group of $S$. For all groups except $\mathsf{B}_2(2)'$, the exceptional part of the Schur multiplier is a $p$-group \cite[Table 24.3]{malletesterman2011}. Therefore, by Lemma \ref{lempext} and \cite[Theorem 7.3]{ruhstorfer2017navarro}, the inductive McKay--Navarro condition is satisfied in defining characteristic for $S$ being one of the groups $\PSL_3(2)$, $\PSL_3(4)$, $\PSL_2(4)$, $\PSL_4(2)$, $\PSU_4(2)$, $\PSU_4(3)$, $^2\mathsf{E}_6(2)$, $\PSL_2(9)$, $\OO_7(3)$, $\PSU_6(2), \mathsf{D}_4(2)=\SO_8^+(2)$.
The remaining groups are $\mathsf{G}_2(2)'$, $^2\mathsf{G}_2(3)'$, $\Sp_6(2)$, $\mathsf{G}_2(3)$, $^2\mathsf{B}_2(8)$, $\mathsf{F}_4(2)$, $^2\mathsf{F}_4(2)'$, and $\mathsf{B}_2(2)'$ \cite[Table 24.3]{malletesterman2011}. 

For $S \in \{ {\Sp}_6(2),$ $\mathsf{G}_2(3),$ ${}^2\mathsf{B}_2(8),$ $ \mathsf{G}_2(2)',$ ${}^2\mathsf{G}_2(3)'\}$, $\Out(G)=\Out(S)$ is cyclic. Let $R \in \Syl_p(G)$, $N=N_G(R)$, and $\Gamma=\Aut(G)_R $. We can explicitly compute the actions of $\Gamma \times \mathcal{H}$ on ${\Irr}_{p'}(G)$ and on ${\Irr}_{p'}(N)$ and see that there exists a $\Gamma \times \mathcal{H}$-equivariant bijection $\Omega$ between the sets. 
We find an outer automorphism $\gamma \in \Aut(G)$ such that $\langle \gamma \rangle \cong \Out(G)$ and can compute the character tables of $G \rtimes \langle \gamma \rangle_\psi$ and $N \rtimes \langle \gamma \rangle_\psi$ for every $\psi \in \Irr_{p'}(G)$. If $G$ is not $\mathsf{G}_2(3)$, we always find $(\Gamma \times \mathcal{H})_\psi$-invariant extensions of $\psi$ and $\Omega(\psi)$. These characters can be canonically extended to the inner automorphisms in $\Gamma$ such that they are still $(\Gamma \times \mathcal{H})_\psi$-invariant. Thus, the inductive McKay--Navarro condition is true for these groups.

If $G=\mathsf{G}_2(3)$, we can construct suitable $(\Gamma \times \mathcal{H})_\psi$-invariant extensions for all $\psi \in \Irr_{3'}(G) \cup \Irr_{3'}(N)$ except for $\chi_2 \in \Irr_{3'}(G)$ of degree $14$ as in \cite[p.60]{ATLAS} and a character $\varphi \in \Irr_{3'}(N)$ of degree $2$. 
The actions of both $\mathcal{H}$ and $\Gamma$ are trivial on these characters; thus we can assume that $\Omega$ maps the characters onto another.
We can now compute the character tables of $G \rtimes \langle \gamma \rangle$ and $N \rtimes \langle \gamma \rangle$. 
Since $N \rtimes \Gamma$ acts trivially on all extensions, we see for both  $\chi_2$ and $\varphi$ 
$$\mu_{y\sigma}= \left\{ \begin{array}{ll}
\chi & \text{ if } {\zeta_3}^\sigma= {\zeta_3}^2, \\ 
1 & \text{ if } {\zeta_3}^\sigma= {\zeta_3}, \\ 
\end{array} \right.$$ 
with $\zeta_3$ a primitive third root of unity and $(y, \sigma) \in (N \rtimes \Gamma) \times \mathcal{H}$. Here, $\chi$ is the character of $G \rtimes \langle \gamma \rangle$ or $N \rtimes \langle \gamma \rangle$, respectively, given by inflation from the non-trivial character of $C_2$. Since we can extend the characters canonically to the inner automorphisms in $\Gamma$, Condition \ref{imck} is satisfied.

\smallskip
For $S=\mathsf{F}_4(2)$, the Schur multiplier has order $2$ and we can therefore consider $G=S$. We see that $R \in \Syl_2(G)$ is self-normalizing and that the character values of the linear characters of $R$ are integers.  
The values of the characters in ${\Irr}_{2'}(G)$ are given in \cite{ATLAS} and we see that $\mathcal{H}$ acts trivially on both ${\Irr}_{2'}(G)$ and ${\Irr}_{2'}(R)$. The outer automorphism group of $G$ is generated by a graph automorphism of order $2$ stabilizing $R$ \cite[Section 2.5]{gorensteinlyonssolomon1998}. 
We can compute the actions of $\gamma$ on the conjugacy classes of $G$ and $N$ and we see that the actions of $\gamma$ on  ${\Irr}_{2'}(G)$ and  ${\Irr}_{2'}(R)$ are permutation isomorphic. 
Thus, for $\psi \in {\Irr}_{2'}(G)$ we now have either $\Gamma_{\psi}=\Inn(R)$ or we can read off the character values of $G \rtimes \langle \gamma \rangle$ from the character table of the split extension $F_4(2):2$ given in \cite{ATLAS}. As before we conclude that the group satisfies Condition \ref{imck}.

\smallskip
We now consider $S=\mathsf{B}_2(2)'\cong \PSL_2(9)$ which has cyclic Schur multiplier of order $6$. 
For $p=3$, this was already treated above; thus let $p=2$. Let $G$ be the $3$-cover of $S$, $R$ a Sylow $2$-subgroup of $G$, and $N=N_G(R)$. 
Note first $\Gamma=\langle \Inn(N), \gamma_1, \gamma_2 \rangle$ with $\gamma_1, \gamma_2 \in \Aut(G)_R$ of order $2$ but there is no subgroup of $\Gamma$ that is isomorphic to the outer automorphism group. Thus, we choose a subgroup $\Gamma' \subseteq \Gamma$ of order $8$ containing representatives of all outer automorphisms of $G$. We can explicitly compute the actions of $\Gamma'$ and $\mathcal{H}$ on ${\Irr}_{2'}(G)$ and ${\Irr}_{2'}(N)$ and construct a $\Gamma \times \mathcal{H}$-equivariant bijection. 
For all $\psi \in {\Irr}_{2'}(G)$, we can apply Lemma \ref{lemmaextensiontrivial} and find $(\Gamma \times \mathcal{H})_\psi$-invariant extensions that are trivial on the centralizer of $N$ in $G \rtimes \Gamma'_\psi$. The local characters are linear and can be extended trivially by setting $\hat{\chi}(x, \gamma)=\chi(x)$ for all $(x,\gamma) \in N \rtimes \Gamma'$. These extensions are $(\Gamma \times \mathcal{H})_\chi$-invariant and we have $\hat{\chi}(g, \gamma)=\chi(g)=1$ for all inner automorphisms $\gamma \in \Gamma'$ induced by $g \in N$. It follows that Condition \ref{imck} holds.

\smallskip

For the Tits group $S=G={}^2\mathsf{F}_4(2)'$, the Schur multiplier is trivial and the extension to the outer automorphism group of order $2$ is again not split. However, we find a representative $\gamma \in \Gamma$ of the outer automorphism with $\ord(\gamma)=4$.
As before, we can construct the $2'$-characters of $G$ and $R$ and see that $\gamma$ acts trivially on them. If $\sigma \in \mathcal{H}$ fixes the fourth root of unity $\zeta_4$,  it fixes all considered characters; otherwise, it acts by interchanging two character pairs of $G$ and $R$ each. Thus, we can find a $\Gamma \times \mathcal{H}$-equivariant bijection.
We can apply Lemma \ref{lemmaextensiontrivial} to get suitable character extensions for the global characters.
The linear local characters can be extended to $\hat{\chi}$ as described before. 
Let $g \in R$ be the element of order $2$ inducing the inner automorphism $\gamma^2$. By explicit computations, we see that we have $\chi(g)=1$ if $\mathcal{H}_\chi=\mathcal{H}$ and $\chi(g)=-1$ else. In the latter case, let $\tau \in \Irr(R \rtimes \langle \gamma \rangle)$ be inflated from a character of $\langle \gamma \rangle$ of order $4$. Then, $(\tau \hat{\chi}) (g,\gamma^2)=1$ and $\tau \hat{\chi}$ is fixed under $\sigma \in \mathcal{H}_\chi =\mathcal{H}_{\zeta_4}$. Thus, we have constructed extensions of the local and global characters such that Condition \ref{imck} is satisfied.
\end{proof}

\newcommand{\etalchar}[1]{$^{#1}$}
\providecommand{\bysame}{\leavevmode\hbox to3em{\hrulefill}\thinspace}
\providecommand{\MR}{\relax\ifhmode\unskip\space\fi MR }
\providecommand{\MRhref}[2]{%
  \href{http://www.ams.org/mathscinet-getitem?mr=#1}{#2}
}
\providecommand{\href}[2]{#2}

\end{document}